\documentclass[12pt]{amsart}
\usepackage{mathrsfs}
\usepackage{amsfonts}
\usepackage{amssymb}
\usepackage{amsfonts, amscd, amsmath, mathrsfs, amssymb, amsthm, amsxtra, bbding, epsfig, graphicx, latexsym, url, mathbbol, bbold}
\usepackage[papersize={9in,10.3in},textwidth=15cm,textheight=20.5cm,centering]{geometry}
\usepackage{enumerate}

\usepackage{xcolor}
\definecolor{cite}{rgb}{0.00,0.60,1.00}
\definecolor{url}{rgb}{1.00,0.10,0.80}
\definecolor{link}{rgb}{0.00,0.00,1.00}
\usepackage[colorlinks,linkcolor=link,urlcolor=url,citecolor=cite,pagebackref,breaklinks]{hyperref}

\hypersetup{
pdfstartpage=1,
pdfstartview=FitH}

\DeclareFontFamily{U}{mathx}{\hyphenchar\font45}
\DeclareFontShape{U}{mathx}{m}{n}{
      <5> <6> <7> <8> <9> <10>
      <10.95> <12> <14.4> <17.28> <20.74> <24.88>
      mathx10
      }{}
\DeclareSymbolFont{mathx}{U}{mathx}{m}{n}
\DeclareMathAccent{\widecheck}{\mathalpha}{mathx}{"71}

 \usepackage{caption} 
\numberwithin{equation}{section}

\allowdisplaybreaks

\newtheorem{theorem}{Theorem}[section]
\newtheorem{lemma}{Lemma}[section]

\newtheorem{definition}{Definition}[section]
\newtheorem{proposition}{Proposition}[section]

\newtheorem{corollary}{Corollary}[section]

\makeatletter
\newcounter{roem}
\renewcommand{\theroem}{\Roman{roem}}

\newcommand{\c@org@eq}{}
\let\c@org@eq\c@equation
\newcommand{\org@theeq}{}
\let\org@theeq\theequation

\newcommand{\setroem}{
\let\c@equation\c@roem
 \let\theequation\theroem}

\newcommand{\setarab}{
\let\c@equation\c@org@eq
\let\theequation\org@theeq}
\makeatother

\newtheorem*{claim*}{Claim}

\theoremstyle{remark}
\newtheorem{remark}{\bf Remark}[section]

\newcommand{\ud}{\mathrm{d}}
\newcommand{\ue}{\mathrm{e}}

\newcommand{\Gal}{\mathrm{Gal}}

\newcommand{\tr}{\mathrm{tr}}
\newcommand{\sym}{\mathrm{sym}}
\newcommand{\Sym}{\mathrm{Sym}}
\newcommand{\rank}{\mathrm{rank}}
\newcommand{\Swan}{\mathrm{Swan}}
\newcommand{\kl}{\mathrm{Kl}}
\newcommand{\Frob}{\mathrm{Frob}}

\newcommand{\Aut}{\mathrm{Aut}}

\newcommand{\SL}{\mathrm{SL}}
\newcommand{\Sp}{\mathrm{Sp}}
\newcommand{\PGL}{\mathrm{PGL}}

\DeclareMathOperator{\Mod}{mod}

\renewcommand{\bmod}[1]{\,(\Mod{ #1})}

\newcommand{\bb}{\mathbf{b}}

\newcommand{\bk}{\mathbf{k}}
\newcommand{\bn}{\mathbf{n}}
\newcommand{\bx}{\mathbf{x}}

\newcommand{\bA}{\mathbf{A}}
\newcommand{\bC}{\mathbf{C}}
\newcommand{\bF}{\mathbf{F}}

\newcommand{\bP}{\mathbf{P}}
\newcommand{\bQ}{\mathbf{Q}}
\newcommand{\bR}{\mathbf{R}}
\newcommand{\bZ}{\mathbf{Z}}

\newcommand{\cA}{\mathcal{A}}
\newcommand{\cB}{\mathcal{B}}

\newcommand{\cE}{\mathcal{E}}
\newcommand{\cF}{\mathcal{F}}

\newcommand{\cK}{\mathcal{K}}
\newcommand{\cL}{\mathcal{L}}
\newcommand{\cM}{\mathcal{M}}
\newcommand{\cN}{\mathcal{N}}
\newcommand{\cP}{\mathcal{P}}

\newcommand{\cV}{\mathcal{V}}
\newcommand{\cW}{\mathcal{W}}

\newcommand{\fc}{\mathfrak{c}}

\newcommand{\fA}{\mathfrak{A}}

\newcommand{\fS}{\mathfrak{S}}

\def\leq{\leqslant}

\def\geq{\geqslant}

\usepackage{graphicx}
\usepackage{tikz}

\usepackage{lipsum}                     % Dummytext
\usepackage{xargs}                      % Use more than one optional parameter in a new commands
\usepackage[colorinlistoftodos,prependcaption,textsize=tiny]{todonotes}
\newcommandx{\unsure}[2][1=]{\todo[linecolor=red,backgroundcolor=red!25,bordercolor=red,#1]{#2}}
\newcommandx{\change}[2][1=]{\todo[linecolor=blue,backgroundcolor=blue!25,bordercolor=blue,#1]{#2}}
\newcommandx{\info}[2][1=]{\todo[linecolor=OliveGreen,backgroundcolor=OliveGreen!25,bordercolor=OliveGreen,#1]{#2}}
\newcommandx{\improvement}[2][1=]{\todo[linecolor=Plum,backgroundcolor=Plum!25,bordercolor=Plum,#1]{#2}}
\newcommandx{\thiswillnotshow}[2][1=]{\todo[disable,#1]{#2}}

\newcommandx{\toreferee}[2][1=]{\todo[linecolor=blue,backgroundcolor=blue!25,bordercolor=blue,#1]{#2}}

\begin{document}

\title[Bilinear forms with trace functions over arbitrary sets]{Bilinear forms with trace functions over arbitrary sets, and applications to Sato--Tate}

\author{Ping Xi}

\address{School of Mathematics and Statistics, Xi'an Jiaotong University, Xi'an 710049, P. R. China}

\email{ping.xi@xjtu.edu.cn}



\begin{abstract} 
We prove non-trivial upper bounds for general bilinear forms with trace functions of bountiful sheaves, where the supports of two variables can be arbitrary subsets in $\bF_p$ of suitable sizes. This essentially recovers the P\'olya--Vinogradov range, and also applies to symmetric powers of Kloosterman sums and Frobenius traces of elliptic curves.
In the case of hyper-Kloosterman sums, we can beat the P\'olya--Vinogradov barrier by combining additive combinatorics with a deep result of Kowalski, Michel and Sawin on sum-products of Kloosterman sheaves.
Two Sato--Tate distributions of Kloosterman sums and Frobenius traces of elliptic curves in sparse families are also concluded.
\end{abstract}

\dedicatory{{\it \small Professor Jingrun CHEN in memoriam}}

\maketitle

\setcounter{tocdepth}{1}

\section{Introduction and main results}\label{sec:Introduction}

\subsection{Backgrounds}
Let $\boldsymbol\alpha=(\alpha_m)$ and $\boldsymbol\beta=(\beta_n)$ be arbitrary complex coefficients with finite supports. Given a target matrix $\phi,$
one usually needs to bound the bilinear form
\begin{align*}
\sum_m\sum_n\alpha_m\beta_n\phi(m,n)
\end{align*}
in terms of the $\ell_2$-norms $\|\boldsymbol\alpha\|_2$ and $\|\boldsymbol\beta\|_2$, and this
occurs quite frequently in many practical applications to number theory, especially in the distribution of primes. 
The idea of bilinear forms in number theory can date back to Vinogradov \cite{Vi37} in his resolution of ternary Goldbach problem for large odd numbers, and he concentrated in the case $\phi(m,n)=\ue(\theta mn)$ with $\theta\in\bR/\bZ,$ where $\ue(z)=\exp(2\pi iz).$ It is also worthwhile to mention the far-reaching work of Chen \cite{Che75} on this occasion, who successfully employed the bilinear structure of remainder terms in linear sieves applying to capture almost primes in short intervals.

This paper mainly concerns some $\phi$ rooted in arithmetic geometry, among which we would like to mention (hyper-) Kloosterman sums, Frobenius traces of elliptic curves at primes, etc.
The resultant estimates for such bilinear forms with $\phi$ can guarantee equidistributions
of Kloosterman sums and Frobenius traces of elliptic curves in families.

Throughout this paper, we take $p$ to be a prime number, and denote by $\bF_p$ the finite field with $p$ elements. We may also identify $\bF_p$ with $\{0,1,\cdots,p-1\}$, so that an interval in $\bF_p$ means a set consisting of some consecutive integers between $0$ and $p-1$.
Let $\cM,\cN\subseteq\bF_p$ be two arbitrary subsets and consider a function $K:\bF_p\rightarrow\bC$. Assume $\boldsymbol\alpha=(\alpha_m)$ and $\boldsymbol\beta=(\beta_n)$ are arbitrary coefficients with supports in $\cM,\cN$, respectively.
We now form the bilinear form
\begin{align}\label{eq:bilinearform}
\cB(\boldsymbol\alpha,\boldsymbol\beta;K)=\sum_{m\in\cM}\sum_{n\in\cN}\alpha_m\beta_nK(mn).
\end{align}
As a direct application of Cauchy's inequality, we have the trivial estimate
\begin{align*}
|\cB(\boldsymbol\alpha,\boldsymbol\beta;K)|\leqslant \|K\|_\infty  \|\boldsymbol\alpha\|_2\|\boldsymbol\beta\|_2(|\cM||\cN|)^{\frac{1}{2}},
\end{align*}
and it turns out to be very crucial in applications to beat the above trivial bound with $|\cM|,|\cN|$ as small as possible.
If $\cM$ and $\cN$ are both intervals in $\bF_p,$ there is a considerable list of estimates for 
$\cB(\boldsymbol\alpha,\boldsymbol\beta;K)$ with various choices of $K$. Here are some selected instances.

\begin{itemize}
\item A very early example with $K(x)=\chi(x+1)$ was considered in depth by Vinogradov \cite{Vi38} and Karatsuba \cite{Ka70}, where $\chi$ is a non-trivial multiplicative character mod $p$. Their motivation is to capture cancellations among the values of $\chi$ at shifted prime arguments.

\item In general, the P\'olya--Vinogradov method allows one to obtain a non-trivial upper bound for $\cB(\boldsymbol\alpha,\boldsymbol\beta;K)$ as long as $|\cM|>p^{\frac{1}{2}+\varepsilon}$, $|\cN|>p^\varepsilon$, provided that
$K$ does not correlate with itself by multiplicative shifts. Fouvry and Michel \cite{FM98} treated the case $K(x)=\ue(f(x)/p)$ with $f$ being a general rational function over $\bF_p$ which is different from polynomials of degrees $0$ and $1$.
In the case of symmetric powers of Kloosterman sums, one may refer to a classical result of Michel \cite{Mi95b}, and much more general trace functions have been considered by Fouvry, Kowalski and Michel \cite{FKM14}.

\item It is highly desirable and quite challenging to beat the P\'olya--Vinogradov barrier, say $\cM$ and $\cN$ are intervals of length around $p^{\frac{1}{2}}.$ As a far-reaching breakthrough, Kowalski, Michel and Sawin \cite{KMS17} succeeded in doing so if $K$ is given by hyper-Kloosterman sums of any fixed rank. They employed the ``shift by $ab$'' trick essentially due to Vinogradov, Karatsuba, Friedlander and Iwaniec, and the problem reduces to proving the geometric irreducibility of some $\ell$-adic sheaves produced by sums of products of many hyper-Kloosterman sums. This was later generalized to situations of more general sums with multiplicative twists; see \cite{KMS20} for details.

\item It is very worthwhile to mention that Shkredov \cite{Sh21} studied the case when $K$ is given by classical Kloosterman sums over $\bF_p$, and employed the growth in $\SL_2(\bF_p)$, estimates for additive energies and tools from incidence geometry. His argument relies heavily on the shape of Kloosterman sums, and gave a purely combinatorial treatment to the involved bilinear forms.
There are also some earlier results by Shparlinski \cite{Sh19} by purely elementary tools, also relying heavily on the shape of Kloosterman sums.
\end{itemize}

\subsection{Main results}
The task of this paper is to bound $\cB(\boldsymbol\alpha,\boldsymbol\beta;K)$ from above when $\cM$ and $\cN$ are general subsets of $\bF_p$. However, if neither of $\cM,\cN$
is contained in suitable intervals, there is no hope to transform incomplete sums to complete sums directly by Fourier analysis (as in the P\'olya--Vinogradov method). By raising powers in the application of H\"older's inequality, we are able to prove non-trivial bounds for $\cB(\boldsymbol\alpha,\boldsymbol\beta;K)$ for some general $K$ (coming from bountiful sheaves as defined in Definition \ref{def:bountifulsheaf}), so that we may also recover the essential range of P\'olya--Vinogradov.

\begin{theorem}\label{thm:arbitraryKMN}
Let $k$ be a positive integer, and $K:\bF_p\rightarrow\bC$ the trace function of some sheaf $\cF$ on $\bA^1_{\bF_p}$ which satisfies one of the following conditions:

$(1)$ $\cF$ is bountiful in the sense of Definition $\ref{def:bountifulsheaf}.$

$(2)$ $K$ is the $k$-th symmetric power of the Kloosterman sum $\kl_2$ given by $\eqref{eq:Kl-tracefunction}.$

$(3)$ $K$ is the $k$-th symmetric power of the Frobenius trace of elliptic curves given by $\eqref{eq:elliptic-tracefunction}.$

Then we have
\begin{align*}
\cB(\boldsymbol\alpha,\boldsymbol\beta;K)
\ll \|\boldsymbol\alpha\|_{\frac{2r}{2r-1}}(p^{\frac{1}{2r}}|\cN|^{\frac{1}{2}-\frac{1}{2r}}\|\boldsymbol\beta\|_{2r}+p^{\frac{1}{4r}}\|\boldsymbol\beta\|_1)
\end{align*}
for any positive integer $r,$ where the implied constant depends only on $r$ and polynomially on the conductor $\fc(\cF).$ Here and henceforth, denote by $\|\cdot\|_q$ the $\ell_q$-norm.
\end{theorem}

One may leave the bountiful sheaf as a black box at this point. To understand Theorem \ref{thm:arbitraryKMN}, we are allowed to take $K$ to be the hyper-Kloosterman sum
\begin{align}\label{eq:hyper-Kl}
\kl_k(a,p)=p^{\frac{1-k}{2}}\mathop{\sum\cdots\sum}_{\substack{x_1,\cdots,x_k\in\bF_p^\times\\ x_1\cdots x_k=a}}\ue\Big(\frac{x_1+\cdots+x_k}{p}\Big).
\end{align}
Note that the work of Deligne \cite{De80} on Riemann Hypothesis for varieties over finite fields guarantees that $|\kl_k(a,p)|\leqslant k$ for all $a\in\bF_p$ and $k\geqslant2.$

If assuming $\boldsymbol\alpha,\boldsymbol\beta$ to be bounded coefficients, the bound in Theorem \ref{thm:arbitraryKMN} then becomes
\begin{align*}
\cB(\boldsymbol\alpha,\boldsymbol\beta;K)
\ll |\cM||\cN|(p^{\frac{1}{2r}}|\cM|^{-\frac{1}{2r}}|\cN|^{-\frac{1}{2}}+p^{\frac{1}{4r}}|\cM|^{-\frac{1}{2r}}),
\end{align*}
which is non-trivial as long as
\begin{align*}
|\cM|>p^{\frac{1}{2}}\log p,\ \ |\cN|>p^\varepsilon
\end{align*}
by taking $r>1/\varepsilon.$ This essentially recovers the range of P\'olya--Vinogradov.

As in many existing works, our next task is then to beat the P\'olya--Vinogradov barrier, at least in some interesting cases. With such intuition in mind, it turns out that we are able to work with hyper-Kloosterman sums as in \cite{KMS17}, but the subsets $\cM,\cN$ here are not necessarily intervals in $\bF_p.$

\begin{theorem}\label{thm:hyperKl}
Let $K=\kl_k(a\cdot,p)$ with $a\in\bF_p^\times$ and $k\geqslant2$. For each positive integer $r\geqslant2$ and arbitrary subsets $\cM,\cN\subseteq\bF_p$ satisfying $|\cN+\cN|\leqslant \lambda|\cN|$ with some $\lambda\geqslant1,$ we have
\begin{align*}
\cB(\boldsymbol\alpha,\boldsymbol\beta;K)
&\ll \|\boldsymbol\alpha\|_\infty\|\boldsymbol\beta\|_2 |\cM||\cN|^{\frac{1}{2}}\Big\{|\cM|^{-\frac{1}{2}}+\Big(\frac{p^{3+\frac{9\lambda}{4r}}\gamma_1\gamma_2}{|\cM|^4|\cN|^3}\Big)^{\frac{1}{8r}}\gamma_3(\log p)^{\frac{1}{2r}}\Big\},
\end{align*}
where $\gamma_1=1+|\cN|^{\frac{3}{2}}p^{-1},\gamma_2=1+|\cN|^{\frac{3}{2}}p^{-1-\frac{9\lambda}{4r}},\gamma_3=1+p^{\frac{3-r}{16r^2}},$ and the implied constant depends only on $(r,\lambda)$ and polynomially on $k$.
\end{theorem}

The condition $|\cN+\cN|\leqslant \lambda|\cN|$ roughly means $\cN$ is not quite far from intervals or arithmetic progressions, and it can be clearly interpreted by the sum-product phenomenon; see Section \ref{sec:thmhyperKl-proof} for details. The proof of Theorem \ref{thm:hyperKl} is inspired by the work of Chang \cite{Cha08} on estimates for the double character sum
\begin{align*}
\sum_m\sum_n\chi(m+n),
\end{align*}
where $m,n$ run over some subsets of $\bF_p.$ The arguments therein combine Burgess' method (essentially the ``shift by $ab$" trick as mentioned above) and a new estimate on the multiplicative energy for subsets due to herself.
In this way, the resultant estimate benefits from the quantitative sum-product phenomenon in finite fields.
We will explain the proof of Theorem \ref{thm:hyperKl} in Section \ref{sec:thmhyperKl-proof} with necessary backgrounds from additive combinatorics.

To see the strength of Theorem \ref{thm:hyperKl}, we may assume $\boldsymbol\alpha,\boldsymbol\beta$ are bounded coefficients for simplicity. By taking $r$ sufficiently large, our bound in Theorem \ref{thm:hyperKl} beats the trivial estimate as long as
\begin{align*}
|\cM|>p^{\varepsilon},\ \ |\cM|^{\frac{4}{3}}|\cN|>p^{1+\varepsilon}.
\end{align*}
This observation allows us to derive the following consequence.

\begin{corollary}\label{coro:hyperKl}
Let $K=\kl_k(a\cdot,p)$ with $a\in\bF_p^\times$ and $k\geqslant2$. Suppose $\|\boldsymbol\alpha\|_\infty,\|\boldsymbol\beta\|_\infty\leqslant1.$ For any fixed $\varepsilon>0$ and $\lambda\geqslant1,$ there exists some positive number $\delta=\delta(\varepsilon,\lambda)$ such that
\begin{align}\label{eq:upperbound-hyperKlcase}
\cB(\boldsymbol\alpha,\boldsymbol\beta;K)
\ll |\cM||\cN|p^{-\delta}
\end{align}
holds for 
\begin{align*}
\cM,\cN\subseteq\bF_p,\ \ |\cN+\cN|\leqslant \lambda|\cN|,\ \ |\cN|\leqslant p^{\frac{2}{3}},\ \ \min\{|\cM|,|\cN|\}>p^{\frac{3}{7}+\varepsilon}.
\end{align*}
\end{corollary}

Note that a very special case of Kowalski, Michel and Sawin \cite[Theorem 4.1]{KMS20} gives the bound \eqref{eq:upperbound-hyperKlcase} for intervals $\cM,\cN$ with
\begin{align*}
\min\{|\cM|,|\cN|\}>p^{\frac{3}{8}+\varepsilon}.
\end{align*}
On the other hand, Bag and Shparlinski \cite{BS23} relaxed one interval to an arbitrary subset of $\bF_p$ under the condition
\begin{align*}
\min\{|\cM|,|\cN|\}>p^{\frac{5}{13}+\varepsilon}.
\end{align*}

\subsection{Sato--Tate distribution of Kloosterman sums and elliptic curves}
Note that $\kl_2$ refers to the classical Kloosterman sum with a normalization, and in Theorem \ref{thm:arbitraryKMN} we can take $K$ to be any symmetric powers of $\kl_2$, which allow us to conclude the following equidistribution of Kloosterman sums.

\begin{theorem}\label{thm:Kl-equidistribution}
Let $p$ be a large prime and $\cM,\cN\subseteq\bF_p$ arbitrary subsets with
\begin{align*}
|\cM|>p^{\frac{1}{2}+\varepsilon},\ \ |\cN|>p^\varepsilon
\end{align*}
for any fixed $\varepsilon>0.$ Then for each $a\in\bF_p^\times,$ the set
\begin{align*}
\{\kl_2(amn):m\in\cM,n\in\cN\}
\end{align*}
becomes equidistributed in $[-2,2]$ with respect to the Sato--Tate measure $\frac{1}{2\pi}\sqrt{4-x^2}\ud x$ as $p\rightarrow+\infty$ over primes.
\end{theorem}

A celebrated result of Katz \cite[Example 13.6]{Ka88} asserts that $\kl_2(am)$ equidistributes in $[-2,2]$, for each fixed $a\in\bF_p^\times,$ with respect to the Sato--Tate measure $\frac{1}{2\pi}\sqrt{4-x^2}\ud x$ as $p\rightarrow+\infty$ over primes. This is usually known as the vertical Sato--Tate law for Kloosterman sums. Michel \cite{Mi95b} was able to prove a similar equidistribution if $m$ only runs over an interval longer than $p^{\frac{1}{2}}(\log p)^2.$ Building on the classical P\'olya--Vinogradov method, one can also succeed in producing equidistributions 
of $\kl_2(amn)$ indexed by $(m,n)\in\cM\times\cN$, as long as $\cM$ and $\cN$ are densely contained in two intervals such that one is longer than $p^{\frac{1}{2}}(\log p)^2$ and the other is longer than $\log p$. See \cite[Corollaire 2.11]{Mi95b} or \cite[Theorem 1.17]{FKM14} for details. The merit of Theorem \ref{thm:Kl-equidistribution} lies in the fact that we are able to deal with arbitrary subsets $\cM,\cN$ with suitably large sizes, agreeing with the 
 P\'olya--Vinogradov barrier.

A similar phenomenon can also be guaranteed in the situation of elliptic curves. To be precise, we consider an elliptic curve $E$ over $\bQ$, and  denote by $E_p$ the reduction of $E$ modulo $p.$ Moreover, denote by $E_p(\bF_p)$ the group of $\bF_p$-rational points of $E_p$. 
The Frobenius trace is defined as usual by 
\begin{align*}
a_p(E)= p+1-|E_p(\bF_p)|.
\end{align*}
The celebrated Hasse bound asserts that $|a_p(E)|\leqslant 2\sqrt{p}.$ The distribution of $a_p(E)$, while $p$ or $E$ runs over suitable families, has received considerable attentions. The horizontal distribution, as $p$ runs over all good primes, was conjectured independently by Sato and Tate in 1960's, and the non-CM case has been settled in a series of papers by Clozel, Harris and Taylor \cite{CHT08}, Taylor \cite{Ta08} and Harris, Shepherd-Barron and Taylor \cite{HST10}.
The vertical analogue seems much easier than the horizontal one, and it was proven by Birch \cite{Bir68} that
$a_p(E)/\sqrt{p}$ equidistributes in $[-2,2]$ with respect to the Sato--Tate measure $\frac{1}{2\pi}\sqrt{4-x^2}\ud x$
as long as $E$ runs over a suitably dense family. The situation becomes quite difficult if $E$ is parametrized in a relatively sparse family. To be precise, we consider the Weierstrass family of elliptic curves
\begin{align*}
E(t):\ y^2=x^3+a(t)x+b(t),
\end{align*}
where and we henceforth assume $a(t),b(t)\in\bZ[t]$ to be polynomials such that
\begin{align*}
\varDelta(t):=-16(4a(t)^3+27b(t)^2)\neq0,\ \ \ j(t):=\frac{-1728(4a(t))^3}{\varDelta(t)}\not\in\bQ.
\end{align*}

\begin{theorem}\label{thm:ellipticcurve-equidistribution}
Let $p$ be a large prime. Let $E(t)/\bQ$ be an elliptic curve as above and denote by $a_p(t)=a_p(E(t))$ the Frobenius trace at $p.$  For each $a\in\bF_p^\times,$ and for arbitrary subsets $\cM,\cN\subseteq\bF_p$ satisfying
\begin{align*}
|\cM|>p^{\frac{1}{2}+\varepsilon},\ \ |\cN|>p^\varepsilon
\end{align*}
with any $\varepsilon>0,$ the set
\begin{align*}
\{a_p(mn)/\sqrt{p}:m\in\cM,n\in\cN,\Delta(mn)\neq0\}
\end{align*}
becomes equidistributed in $[-2,2]$ with respect to the Sato--Tate measure $\frac{1}{2\pi}\sqrt{4-x^2}\ud x$ as $p\rightarrow+\infty$ over primes.
\end{theorem}

There is a similar history on the vertical distribution of Frobenius traces of elliptic curves, and 
the first result on the equidistribution of $a_p(m)/\sqrt{p}$, as $m$ runs over $\bF_p^\times$ for fixed large prime $p$, is due to Birch \cite{Bir68} as mentioned above. One can also obtain similar equidistributions when $m$ runs over intervals or with some bilinear structures thanks to the work in \cite{Mi95a,Mi95b,FKM14}.
The equidistribution in Theorem \ref{thm:ellipticcurve-equidistribution} has been achieved by de la Bret\`eche, Sha, Shparlinski and Voloch \cite{BSSV18} in the situation
\begin{align*}
|\cM||\cN|>p
\end{align*}
with an explicit rate of convergence. A special consequence of Theorem \ref{thm:arbitraryKMN} allows us to find an alternative restriction on the sizes of $\cM,\cN$ such that the desired equidistribution holds.

The idea in proving Theorem \ref{thm:arbitraryKMN} should admit generalizations to a large class of bilinear forms of the shape \eqref{eq:bilinearform}. It came to the author while studying equidistributions of Jacobi sums in \cite{Xi18}, for which we transformed the problem to twisted moments of Gauss sums, and we would come back to this subject in a forthcoming paper. On the other hand, the proof of Theorem \ref{thm:hyperKl} combines additive combinatorics with a deep result of Kowalski, Michel and Sawin \cite{KMS20} on the {\it sum-product} of Kloosterman sheaves; see Lemma \ref{lm:KMS-hyperKl} for details. In fact, their estimate was originally formulated in terms of general hyper-Kloosterman sums with multiplicative twists:
\begin{align*}
\kl_k(a,\boldsymbol\chi;p)=p^{\frac{1-k}{2}}\mathop{\sum\cdots\sum}_{\substack{x_1,\cdots,x_k\in\bF_p^\times\\ x_1\cdots x_k=a}}\chi_1(x_1)\cdots \chi_k(x_k)\ue\Big(\frac{x_1+\cdots+x_k}{p}\Big),
\end{align*}
where $\boldsymbol\chi=(\chi_1,\cdots,\chi_k)$ is a tuple of multiplicative characters over $\bF_p^\times.$
Our bound in Theorem \ref{thm:hyperKl} also generalizes to this framework, but we choose to keep our situation in hyper-Kloosterman sums \eqref{eq:hyper-Kl} avoiding a couple of new concepts.

This note is organized as follows. We introduce some terminology on trace functions and sheaves, in particular bountiful sheaves, in Section \ref{sec:tracefunctionssheaves}.
The proofs of Theorem \ref{thm:arbitraryKMN} and \ref{thm:hyperKl} will be given in Sections \ref{sec:thmarbitraryKMN-proof} and \ref{sec:thmhyperKl-proof}, respectively.
The equidistributions in Theorems \ref{thm:Kl-equidistribution} and \ref{thm:ellipticcurve-equidistribution} will be proved in the last section.

\smallskip

\subsection*{Acknowledgements} 
I am very grateful to Bryce Kerr for letting me know the work of Shkredov \cite{Sh18} which yields an improved version of Theorem \ref{thm:hyperKl}, and to Igor Shparlinski for pointing out the reference \cite{BS23}. 
I also thank the referees for valuable comments and suggestions.

It is my great honour to be invited to acknowledge the 50th anniversary of the detailed proof of Chen's celebrated theorem on the Goldbach problem and twin prime conjecture. I am very lucky that my academic career has been guided by Chen's mathematics and spirits.

This work is supported in part by NSFC (No. 12025106, No. 11971370).

\smallskip

\section{Trace functions and sheaves}
\label{sec:tracefunctionssheaves}

In this section, we introduce the terminology on trace functions of $\ell$-adic sheaves on $\bA_{\bF_p}^1$ 
following the manner of Fouvry, Kowalski and Michel \cite{FKM14, FKM15}.

\subsection{Trace functions}
Let $p$ be a prime and $\ell\neq  p$ an auxiliary prime, and fix an isomorphism $\iota : \overline{\bQ}_\ell\rightarrow\bC$. The functions $K(x)$ modulo $p$ that we
consider are the trace functions of suitable constructible sheaves on $\bA^1_{\bF_p}$
evaluated at $x\in\bF_p$. To
be precise, we will consider middle-extension sheaves on $\bP^1_{\bF_p}$ 
and we refer to the following definition after Katz \cite[Section 7.3.7]{Ka88}.

\begin{definition}[Trace functions]\label{def:tracefunction}
Let $\cF$ be an $\ell$-adic middle-extension sheaf pure of weight zero, 
which is lisse on an open set $U$. The trace function associated to $\cF$ is defined by
\begin{align*}
K:x\in\bF_p\mapsto\iota(\tr(\Frob_x\mid V_\cF)),
\end{align*}
where $\Frob_x$ denotes the geometric Frobenius at $x\in\bF_p,$ and $V_\cF$ is a finite dimensional $\overline{\bQ}_\ell$-vector space, which is corresponding to a continuous finite-dimensional Galois representation and unramified at every closed point $x$ of $U.$
\end{definition}

We need an invariant to measure the geometric complexity of a trace function, which can be given by some numerical invariants of the underlying sheaf.
\begin{definition}[Conductor] \label{def:conductor} 
For an $\ell$-adic middle-extension sheaf $\cF$ on $\bP^1_{\bF_p}$ of rank $\rank(\cF)$,
we define the $($analytic$)$ conductor of $\cF$ to be
\begin{align*}  
\fc(\cF) := \rank(\cF) + \sum_{x\in S(\cF)} (1+\Swan_x(\cF)),
\end{align*}
where $S(\cF)\subset\bP^1(\overline{\bF}_p)$ denotes the $($finite$)$ set of singularities of $\cF$, 
and $\Swan_x(\cF)$ $(\geqslant 0)$ denotes the Swan conductor of $\cF$ at $x$ $($see {\rm \cite{Ka80}}$).$
\end{definition}

\subsection{Bountiful sheaves}
A large body of this paper concerns the following special sheaves, and the definition is borrowed directly from Fouvry, Kowalski and Michel \cite{FKM15}.

\begin{definition}[Bountiful sheaves]\label{def:bountifulsheaf}
An $\ell$-adic sheaf $\cF$ on $\bA^1_{\bF_p}$, which is middle extension and pointwise pure of
weight $0$, is said to be bountiful if the following conditions hold:
\begin{itemize}
\item The rank of $\cF$ is at least $2$;
\item The arithmetic monodromy group of $\cF$ is equal to the geometric monodromy group, and is equal to either $\SL_r$ or $\Sp_{r}$;
\item The projective automorphism group
\begin{align*}
\Aut_0(\cF)=\{\gamma\in\PGL_2(\bF_p)\,\mid\,
\gamma^*\cF\simeq \cF\otimes \cL\text{ for some rank
  $1$ sheaf }\cL\}
\end{align*}
of $\cF$ is trivial.
\end{itemize}
We will say $\cF$ is of $\SL$-type or $\Sp$-type accordingly.
\end{definition}

There are many instances of $\SL_r$-type and $\Sp_r$-type bountiful sheaves. Recall the definition \eqref{eq:hyper-Kl} of hyper-Kloosterman sums. As computed by Katz \cite[Theorem 11.1]{Ka88}, for each $k\geqslant2$ and $p>2$, $\kl_k(\cdot,p)$ is the trace function of a bountiful sheaf $\cK l_k$ on $\bA_{\bF_p}^1$, which is of $\SL_k$ type if $k$ is odd, and of $\Sp_k$ type if $k$ is even.
More details and examples can be found in \cite{FKM15}.

\begin{definition}[Normal tuples]\label{def:normal}
Let $p$ be a prime, $k$ a positive integer, $\boldsymbol\gamma$ a $k$-tuple of $\PGL_2(\bF_p)$ and $\boldsymbol\sigma$ a $k$-tuple of $\Gal(\bC/\bR)=\{1,c\}$, where $c$ is complex conjugation.

$(1)$ We say that $\boldsymbol\gamma$ is normal if there exists some
$\gamma\in\PGL_2(\bF_p)$ such that
$$
|\{1\leq i\leq k\,\mid\, \gamma_i=\gamma\}|\equiv1\bmod2.
$$

$(2)$ If $r\geq 3$ is an integer, we say that
$(\boldsymbol\gamma,\boldsymbol\sigma)$ is $r$-normal if there exists
some $\gamma\in\PGL_2(\bF_p)$ such that
$$
|\{1\leq i\leq k\,\mid\, \gamma_i=\gamma\}|\geq 1
$$
and
$$
|\{1\leq i\leq k\,\mid\, \gamma_i=\gamma\text{ and } \sigma_i=1 \}|\not\equiv
|\{1\leq i\leq k\,\mid\, \gamma_i=\gamma\text{ and } \sigma_i\not=1
\}|\bmod r.
$$
\end{definition}

Given an $\ell$-adic sheaf $\cF$ with the trace function $K$ and with conductor $\fc$, we consider the sum of product
\begin{align}\label{eq:sumofproduct}
\fS(K;\boldsymbol\gamma,\boldsymbol\sigma)
:=\sum_{x\in\bF_p}K(\gamma_1\cdot x)^{\sigma_1}\cdots K(\gamma_k\cdot x)^{\sigma_k},
\end{align} 
where $k\geqslant1$, $\boldsymbol\gamma$ and $\boldsymbol\sigma$ are $k$-tuples of elements of $\PGL_2(\bF_p)$ and $\Gal(\bC/\bR)$, respectively. For each $\sigma\in\Gal(\bC/\bR),$ we write $K^\sigma=\sigma(K).$
The average $\fS(K;\boldsymbol\gamma,\boldsymbol\sigma)$ appears in many arithmetic problems in analytic number theory, which requires one to capture the square root cancellation as long as $\boldsymbol\gamma$ and $\boldsymbol\sigma$ appear in suitable configurations.

Fouvry, Kowalski and Michel \cite{FKM15} proved the following two estimates for bountiful sheaves $\cF$ in full generality; see Corollary 1.6 and Corollary 1.7 therein.

\begin{proposition}\label{prop:noncorrelation}
Let $p$ be a prime and let $K$ be the trace function modulo $p$ of a
bountiful sheaf $\cF$ with conductor $\fc$. Then for any $k\geqslant1$, and for any $k$-tuple $\boldsymbol\gamma$ of elements of $\PGL_2(\bF_p)$ and $\boldsymbol\sigma$ of $\Gal(\bC/\bR)$, the inequality
\begin{align*}
|\fS(K;\boldsymbol\gamma,\boldsymbol\sigma)|\leqslant C\sqrt{p}
\end{align*}
holds with a constant $C=C(k,\fc)$ depending only on $k,\fc$, provided that one of the following conditions is satisfied:

$(1)$ The sheaf $\cF$ is self-dual (so that $K$ is real-valued) and $\boldsymbol\gamma$ is normal.

$(2)$ The sheaf $\cF$ is of $\SL_r$-type with $p>r\geqslant3$, and
$(\boldsymbol\gamma,\boldsymbol\sigma)$ is $r$-normal.

\noindent The dependence of $C$ on $\fc$ is polynomial.
\end{proposition}

\begin{proposition}\label{prop:correlation}
Let $p$ be a prime and let $K$ be the trace function modulo $p$ of a bountiful sheaf $\cF$ with conductor $\fc$. 
Then for any $k\geqslant1$, and for any $k$-tuple $\boldsymbol\gamma$ of elements of $\PGL_2(\bF_p)$ and $\boldsymbol\sigma$ of $\Gal(\bC/\bR)$, there exists some positive integer $m(\boldsymbol\gamma,\boldsymbol\sigma)\geqslant1$ depending on $k,\fc$, such that the inequality
\begin{align*}
|\fS(K;\boldsymbol\gamma,\boldsymbol\sigma)-m(\boldsymbol\gamma,\boldsymbol\sigma)p|\leqslant C\sqrt{p}
\end{align*}
holds with a constant $C=C(k,\fc)$ depending only on $k,\fc$, provided that one of the following conditions is satisfied:

$(1)$ The sheaf $\cF$ is of $\Sp_{2g}$-type and $\boldsymbol\gamma$ is not
normal.

$(2)$ The sheaf is of $\SL_r$-type with $r\geqslant3$,
$(\boldsymbol\gamma,\boldsymbol\sigma)$ is not $r$-normal, and the pull-back $[x\mapsto-x]^*\cF$ is geometrically isomorphic to the dual of $\cF,$ i.e.,
\begin{align*}
[x\mapsto-x]^*\cF\simeq D(\cF).
\end{align*}

\noindent The dependence of $C$ on $\fc$ is polynomial.
\end{proposition}

Proposition \ref{prop:noncorrelation} establishes an asymptotic formula for suitable $\boldsymbol\gamma,\boldsymbol\sigma$ as $p\rightarrow+\infty$. This result is merely presented here for completeness: we will use the trivial bound $\fS(K;\boldsymbol\gamma,\boldsymbol\sigma)\ll p$ when $\boldsymbol\gamma,\boldsymbol\sigma$ appear in some particular configurations, since we are seeking upper bound, instead of asymptotic formulae, for bilinear forms.

\subsection{Correlations among Kloosterman sums and Frobenius traces}
Following Deligne \cite{De80} and Katz \cite{Ka88}, it is known that
\begin{align}\label{eq:Kl-tracefunction}
a\mapsto-\kl_2(a,p)=-2\cos\theta_p(a),\ \ a\in \bF_p^\times
\end{align}
is the trace function of an $\ell$-adic sheaf $\mathcal{K}l$ on $\mathbf{G}_{m}(\bF_p)=\bF_p^\times$, which is of rank 2 and pure of weight 0. 
Alternatively, we may write
\begin{align*}2\cos\theta_p(a)=\tr(\Frob_a,\mathcal{K}l_2),\ \ a\in \bF_p^\times.\end{align*}
Note that the arithmetic monodromy group is equal to the geometric monodromy group, and both are $\SL_2=\Sp_2$ according to Katz \cite[Theorem 11.1]{Ka88}.
Put
\begin{align*}\sym_k(\theta)=\frac{\sin(k+1)\theta}{\sin\theta},\end{align*}
so that $a\mapsto \sym_k(\theta_p(a)) $
is the trace function of $\Sym^k\mathcal{K}l_2$, the $k$-th symmetric power
of the Kloosterman sheaf $\mathcal{K}l$ (i.e., the composition of the sheaf $\mathcal{K}l$ with the $k$-th symmetric power representation of $\SL_2$).
In this way, the geometric monodromy group of $\Sym^k\mathcal{K}l$ is $\Sym^k(\SL_2).$

The choice $\cF=\Sym^k\mathcal{K}l_2$ of course does not fall into the framework of bountiful sheaves. However, the representations of $\Sym^k(\SL_2)$ are very clearly understood, so that we may also conclude the estimates in Propositions \ref{prop:noncorrelation} and \ref{prop:correlation} in the situation of symmetric powers of Kloosterman sums.

To be precise, we put
\begin{align}\label{eq:correlationsum-Kloosterman}
\fA(\bk,\boldsymbol\gamma;h)
&:=\sum_{x\in\bF_p}\prod_{1\leqslant j\leqslant s}\sym_{k_j}(\theta_p(\gamma_j\cdot x))\ue\Big(\frac{hx}{p}\Big),
\end{align}
where $s$ is a positive integer, $\bk$ is an $s$-tuple of positive integers and $\boldsymbol\gamma$ is an $s$-tuple of elements of $\PGL_2(\bF_p)$.
In fact, Katz \cite[Example 13.6]{Ka88} proved that
\begin{align}\label{eq:Katz}
\left|\sum_{a\in\bF_p^\times}\sym_k(\theta_p(a))\right|\leqslant\frac{1}{2}(k+1)\sqrt{p},
\end{align}
which allows him to conclude the vertical Sato--Tate law from Weyl's criterion. This gives an upper bound for 
$\fA(\bk,\boldsymbol\gamma;h)$ in the case $s=1,\gamma=\mathrm{Id}$ and $h=0$.
A non-trivial estimate for $\fA(\bk,\boldsymbol\gamma;h)$ in generic cases can be regarded as a high-dimensional analogue of \eqref{eq:Katz}, and this follows from the independence of
monodromy groups of Kloosterman sheaves (as well as their symmetric powers).

\begin{proposition}\label{prop:noncorrelation-Kloosterman}
Let $p$ be a prime. For any $s\geqslant1$, and for any $s$-tuple $\bk$ of positive integers and $s$-tuple $\boldsymbol\gamma$ of elements of $\PGL_2(\bF_p),$ the inequality
\begin{align*}
|\fA(\bk,\boldsymbol\gamma;h)|\leqslant c\prod_{1\leqslant j\leqslant s}k_j^2\cdot \sqrt{p}
\end{align*}
holds with a constant $c=c(s)$ depending only on $s,$ provided that one of the following conditions is satisfied:

$(1)$ $h\in\bF_p^\times.$

$(2)$ There exists some $\gamma\in\PGL_2(\bF_p)$ such that
\begin{align*}
\sum_{\substack{1\leqslant j\leqslant s\\\gamma_j=\gamma}}k_j\equiv1\bmod2.
\end{align*}
\end{proposition}

A prototype of Proposition \ref{prop:noncorrelation-Kloosterman} was proven by
Fouvry, Michel, Rivat and S\'{a}rk\"{o}zy \cite{FMRS04} in their investigations of the pseudorandomness of signs of Kloosterman sums, where they assume the coordinates of $\boldsymbol\gamma$ are pairwise distinct and all are given by upper triangular matrices in $\PGL_2(\bF_p)$. The case $\bk=(1,\cdots,1)$, which falls into the framework of Proposition \ref{prop:noncorrelation}, is proven previously by Fouvry, Ganguly, Kowalski and Michel \cite[Proposition 3.2]{FGKM14}. The simpler version of Proposition \ref{prop:noncorrelation-Kloosterman}, when $\bk$ is general and $\boldsymbol\gamma$ is composed by upper triangular matrices in $\PGL_2(\bF_p)$,
was obtained by the author \cite[Lemma 4]{Xi17} using combinatorial properties of Chebyshev polynomials together with \cite[Lemma 2.1]{FMRS04}. The arguments therein also apply to general tuples $\boldsymbol\gamma$ of elements of $\PGL_2(\bF_p)$.

In the proof of Theorem \ref{thm:arbitraryKMN}, we only require the case $h=0$, $\bk=(k,\cdots,k)$ and $\boldsymbol\gamma$ being composed by diagonal matrices.
We formulate Proposition \ref{prop:noncorrelation-Kloosterman} in its current version since it is of independent interests and might be applicable to some other problems.

We now turn to the situation of elliptic curves, which have similar interpretations to those of Kloosterman sums.
Recall the notation in the first section. It is known that
\begin{align}\label{eq:elliptic-tracefunction}
t\mapsto a_p(t)/\sqrt{p}=2\cos\widetilde{\theta}_p(t),\ \ t\in \bF_p, \varDelta(t)\neq0
\end{align}
is the trace function of an $\ell$-adic sheaf $\cE$ on $\bF_p-\{\varDelta(x)=0\}$, which is of rank 2 and pure of weight 0. 
We may alternatively write
\begin{align*}
2\cos\widetilde{\theta}_p(t)=\tr(\Frob_t,\cE),\ \ t\in \{x\in\bF_p:\varDelta(x)\neq0\}.
\end{align*}
According to Deligne \cite[Lemme (3.5.5)]{De80}, the geometric monodromy group is equal to the arithmetic monodromy group, $\SL_2$. Following the discussions on Kloosterman sums as above, we can also speak of symmetric powers of $\cE,$ and the estimate in Proposition \ref{prop:noncorrelation-Kloosterman} also holds with an upper bound
\begin{align*}
c\cdot r_\varDelta^s\cdot \prod_{1\leqslant j\leqslant s}k_j^2\cdot \sqrt{p}
\end{align*}
if replacing 
$\theta_p$ by $\widetilde{\theta}_p$ in \eqref{eq:correlationsum-Kloosterman}, where $r_\varDelta=|\{z\in\bC:\varDelta(z)=0\}|$.

\smallskip

\section{Proof of Theorem \ref{thm:arbitraryKMN}: P\'olya--Vinogradov range}
\label{sec:thmarbitraryKMN-proof}

Let $r$ be a positive integer. By H\"older's inequality, we have
\begin{align}\label{eq:Holder-initial}
|\cB(\boldsymbol\alpha,\boldsymbol\beta;K)|\leqslant \|\boldsymbol\alpha\|_{2r/(2r-1)}\cB^{1/2r},
\end{align}
where
\begin{align*}
\cB=\sum_{m\in\cM}\Big|\sum_{n\in\cN}\beta_nK(mn)\Big|^{2r}.
\end{align*}

We amplify the $m$-sum in $\cB$ to the whole $\bF_p$, so that
\begin{align}\label{eq:B-B(n)}
\cB\leqslant \sum_{m\in\bF_p}\Big|\sum_{n\in\cN}\beta_nK(mn)\Big|^{2r}
=\sum_{\bn\in\cN^{2r}}\beta(\bn)\cB(\bn),
\end{align}
where for $\bn=(n_1,n_2,\cdots,n_{2r})\in\cN^{2r},$
\begin{align*}
\beta(\bn)=\prod_{1\leqslant j\leqslant r}\beta_{n_j}\overline{\beta}_{n_{j+r}},
\end{align*}
and
\begin{align*}
\cB(\bn)=\sum_{x\in\bF_p}\prod_{1\leqslant j\leqslant r}K(n_jx)\overline{K(n_{j+r}x)}.
\end{align*}

For each $j\geqslant1,$ consider
\begin{align*}
\gamma_j=\Big(\begin{matrix}
n_j & 0\\
0 & 1
\end{matrix}\Big),
\end{align*}
so that 
\begin{align*}
\cB(\bn)=\sum_{x\in\bF_p}\prod_{1\leqslant j\leqslant r}K(\gamma_j\cdot x)\overline{K(\gamma_{j+r}\cdot x)}=\fS(K;\boldsymbol\gamma,\boldsymbol\sigma)
\end{align*}
as defined in \eqref{eq:sumofproduct}, where
\begin{align*}
\boldsymbol\gamma=(\gamma_1,\gamma_2,\cdots,\gamma_{2r}),\ \ \boldsymbol\sigma=( \underbrace{1,1,\cdots,1}_{r\text{ copies}},\underbrace{c,c,\cdots,c}_{r\text{ copies}})
\end{align*} 
with $c$ denoting the complex conjugation.
We would like to bound $\cB(\bn)$ by appealing to Proposition \ref{prop:noncorrelation} if 
$\boldsymbol\gamma,\boldsymbol\sigma$ don't appear in certain configurations, and otherwise the trivial bound
\begin{align*}
\cB(\bn)\ll p
\end{align*}
is applied. To do so, we assume the sheaf $\cF$ is bountiful as defined in Definition \ref{def:bountifulsheaf}.

For $\Sp$-type $\cF$, if all its coordinates appear in pairs (possibly with multiples), 
$\boldsymbol\gamma$ is normal according to Definition \ref{def:normal}, in which case we infer 
\begin{align*}
\cB(\bn)\ll \sqrt{p}
\end{align*}
from Proposition \ref{prop:noncorrelation}. Applying the trivial bound for $\cB(\bn)$ in remaining situations, we find
\begin{align*}
\cB
&\ll p\sum_{\substack{A\subseteq\bZ^+\\ \sum_{a\in A}a=r}}\prod_{a\in A}\|\boldsymbol\beta\|_{2a}^{2a}+\sqrt{p}\|\boldsymbol\beta\|_1^{2r},
\end{align*}
where the elements of $A$ are not necessarily distinct.
By H\"older's inequality, we infer
\begin{align*}
\|\boldsymbol\beta\|_{2a}^{2a}\leqslant N^{1-\frac{a}{r}}\|\boldsymbol\beta\|_{2r}^{2a},
\end{align*}
so that
\begin{align*}
\sum_{\substack{A\subseteq\bZ^+\\ \sum_{a\in A}a=r}}\prod_{a\in A}\|\boldsymbol\beta\|_{2a}^{2a}
&\leqslant N^{r-1}\sum_{\substack{A\subseteq\bZ^+\\ \sum_{a\in A}a=r}}\|\boldsymbol\beta\|_{2r}^{2a}
\ll N^{r-1}\|\boldsymbol\beta\|_{2r}^{2r}.
\end{align*}
This yields
\begin{align}\label{eq:B-upperbound}
\cB
&\ll pN^{r-1}\|\boldsymbol\beta\|_{2r}^{2r}+\sqrt{p}\|\boldsymbol\beta\|_1^{2r}.
\end{align}
Theorem \ref{thm:arbitraryKMN}, in the case of $\Sp$-type $\cF$, now follows immediately from \eqref{eq:Holder-initial}, \eqref{eq:B-B(n)} and \eqref{eq:B-upperbound}.
Following similar arguments, we may also establish the inequality \eqref{eq:B-upperbound} for $\SL$-type $\cF$.
The cases of symmetric powers of Kloosterman sums and Frobenius traces of elliptic curves can be treated similarly using Proposition \ref{prop:noncorrelation-Kloosterman} and the subsequent comments.

\smallskip

\section{Proof of Theorem \ref{thm:hyperKl}: Combinatorial arguments}
\label{sec:thmhyperKl-proof}

The proof of Theorem \ref{thm:hyperKl} benefits from quite a lot of progresses in additive combinatorics, in particular the sum-product phenomenon in finite fields. We first introduce some basic concepts and preliminary results.

\subsection{Sumsets and generalized arithmetic progressions}
For subsets $\cA,\cB\subseteq \bZ$ or $\bF_p$, we define the sumset $\cA+\cB$, difference set $\cA-\cB$ and product set $\cA\cdot\cB$ by
\begin{align}\label{eq:sumdifferenceproductset}
\cA*\cB:=\{a*b:a\in  \cA,b\in\cB\},\ \ \ *\in\{+,-,\cdot\}.
\end{align}

Given integers $a_0, a_1, \ldots a_k$ and positive integers $N_1, N_2, \ldots N_k$, the set
\begin{align*}
\cP:=\{a_0+a_1 n_1+a_2 n_2+\cdots+a_d n_d: 0 \leqslant n_j \leqslant N_j-1 \text { for } 1 \leqslant j \leqslant d\}
\end{align*}
is called a {\it generalized arithmetic progression} of dimension $d$ and volume $N_1 N_2 \ldots N_d$. If $|\cP|=N_1 N_2 \ldots N_d$, i.e., all such numbers are distinct, then $\cP$ is said to be {\it proper}. 

As in Chang \cite{Cha08}, the proof of Theorem \ref{thm:hyperKl} requires Freiman's theorem on sumsets, which asserts that any finite subset $\cA\subseteq\bZ$ should be contained in 
a generalized arithmetic progression, provided that $|\cA+\cA|/|\cA|$ is small.
Besides its origin \cite{Fr73}, one may also see many improvements and developments by Ruzsa \cite{Ru94}, Bilu \cite{Bil99}, Chang \cite{Cha02}, Schoen \cite{Sc11} and Sanders \cite{Sa13} for instance. One may also find a detailed history from \cite{Sa13} on the quantitative dependences of dimension and volume of the progression on the doubling constant $|\cA+\cA|/|\cA|$.
The following result is obtained by Chang \cite[Theorem 2]{Cha02}.

\begin{lemma}\label{lm:Chang-progression}
Assume that $\cA\subseteq\bZ$ is a finite set with $|\cA+\cA|\leqslant \lambda |\cA|$ for some $\lambda\geqslant1.$ Then $\cA$ is contained in a proper $d$-dimensional arithmetic progression $\cP$ with
\begin{align*}
d \leqslant\lambda-1, \ \ \log(|\cP|/|\cA|)\leqslant C \lambda^2(\log \lambda)^3
\end{align*}
for some absolute positive constant $C.$
\end{lemma}

We choose to employ the above work of Chang \cite{Cha02} since the dependence of $d$ on $\lambda$ is linear. There are many other works, as mentioned above, by weakening the dependence of $|\cP|/|\cA|$ on $\lambda$ at the cost of a much larger dimension $d$.

As usual we define the \emph{multiplicative energy} of $\cA,\cB:$
\begin{align*}
E(\cA ,\cB):=|\{(a_1,a_2,b_1,b_2)\in  \cA^2 \times \cB^2:~a_1b_1=a_2b_2\}|.
\end{align*}
We write  $E(\cA,\cA)=E(\cA)$ for abbreviation. By orthogonality and Cauchy inequality, one may derive that
\begin{align*}
E(\cA ,\cB)^2\leqslant E(\cA)E(\cB).
\end{align*}
Moreover, define
\begin{align*}
D(\cA):=|\{(a_1,\cdots,a_8)\in  \cA^8:~(a_1-a_2)(a_3-a_4)=(a_5-a_6)(a_7-a_8)\}|.
\end{align*}
The quantity $D(\cA)$ can be interpreted as the number of incidences between points and planes
\begin{align*}
(a_1-a_2) \alpha=(a_5-a_6) \beta
\end{align*}
counting with the weights $|\{x-y=t:x,y\in \cA\}|$ for $t=\alpha,\beta$.
By using incidence theorems in $\bF_p,$ Shkredov \cite[Theorem 32]{Sh18} was able to prove a very strong upper bound for $D(\cA),$ which even offers an asymptotic formula when $\cA$ is relatively large, say $|\cA|>p^{\frac{2}{3}}\log^3p.$
\begin{lemma}\label{lm:Shkredov}
For any subset $\cA\subseteq\bF_p,$ we have
\begin{align*}
D(\cA)=\frac{|\cA|^8}{p}+O(|\cA|^{\frac{13}{2}}\log^2|\cA|).
\end{align*}
\end{lemma}
If $\cA$ is a generalized arithmetic progression, we may give a very good estimate for $E(\cA)$ by virtue of Lemma \ref{lm:Shkredov}.

\begin{lemma}\label{lm:Shkredov-Kerr}
For each generalized arithmetic progression $\cP$ of dimension $d$ in $\bF_p,$ we have
\begin{align*}
E(\cP)\ll 2^{8d}\Big(\frac{|\cP|^4}{p}+|\cP|^{\frac{5}{2}}\log^2|\cP|\Big).
\end{align*}
\end{lemma}

\proof
For each $x\in\cP,$ there are at least $\gg |\cP|$ tuples of $(x_1,x_2)\in(\cP+\cP)^2$ such that $x=x_1-x_2.$
Hence
\begin{align*}
E(\cP)\ll \frac{D(\cP+\cP)}{|\cP|^4}.
\end{align*}
Now the lemma follows directly from Lemma \ref{lm:Shkredov} and $|\cP+\cP|\leqslant 2^d|\cP|$.
\endproof

\begin{remark}
If $\cP$ is a lattice generated by a basis of $\bF_{p^d}/\bF_p$ with $|\cP|<p^{\frac{1}{2}}$, Chang \cite[Proposition 1]{Cha08} was able to give an upper bound
\begin{align*}
E(\cP)\leqslant 5^d|\cP|^{\frac{11}{4}}\log p.
\end{align*}
The arguments therein can also give essentially the same bound if $\cP$ in a generalized arithmetic progression; see (3.2) therein. The work of Shkredov (Lemma \ref{lm:Shkredov}) allows us to derive a better estimate for $E(\cP)$ for all $|\cP|\ll p^{\frac{1}{2}}.$ I am very grateful to Bryce Kerr for letting me know the work of Shkredov \cite{Sh18} and showing me the proof of Lemma \ref{lm:Shkredov-Kerr}. Note that the work of Chang \cite{Cha08} yields a weaker bound in Theorem \ref{thm:hyperKl}, so that produces a more restrictive range $\min\{|\cM|,|\cN|\}>p^{\frac{6}{13}+\varepsilon}$ in Corollary \ref{coro:hyperKl}.
\end{remark}

We are now ready to prove Theorem \ref{thm:hyperKl}. We assume $K(x)=\kl_k(ax,p)$ for $a\in\bF_p^\times$ and $\|\boldsymbol\alpha\|_\infty\leqslant1.$
Let $\cN_1$ be a generalized $d$-dimensional proper arithmetic progression in $\bF_p$ given by
\begin{align*}
\cN_1=\Big\{\sum_{1\leqslant j\leqslant d} x_j \omega_j: 0\leqslant  x_j\leqslant H_j, 1\leqslant j\leqslant d\Big\},
\end{align*}
satisfying $\cN\subseteq a_0+\cN_1$ for some $a_0\in\bF_p,$ and
\begin{align*}
d \leqslant \lambda, \ \ |\cN_1|\leqslant \ue^{C\lambda^2(\log\lambda)^3}|\cN|,
\end{align*}
where $C$ is the positive constant appearing in Lemma \ref{lm:Chang-progression}. 

Let $V$ be a large number to be chosen later, and put
\begin{align*}
I&=[1,V]\cap\bZ,\\
\cN_0&=\Big\{\sum_{1\leqslant j\leqslant d} x_j \omega_j: 0\leqslant  x_j\leqslant V^{-1}H_j, 1\leqslant j\leqslant d\Big\},\\
\cN_2&=a_0+(\cN_1\cup(-\cN_1)),
\end{align*}
so that $\cN_0$ is also a proper progression and
\begin{align}
V^d|\cN_0|\asymp |\cN_2|, \ \ \cN-\cN_0 I \subseteq \cN_2,
\end{align}
where $\cN-\cN_0 I$ is understood as in \eqref{eq:sumdifferenceproductset}.
Note that
\begin{align*}
|\cN| \leqslant|\cN_1|\leqslant \ue^{C\lambda^2(\log\lambda)^3}|\cN| ,\ \ \ |\cN_2|=2|\cN_1|-1.
\end{align*}

By Cauchy, we have
\begin{align*}
|\cB(\boldsymbol\alpha,\boldsymbol\beta;K)|^2
&\leqslant \|\boldsymbol\beta\|^2\sum_{n\in\cN}\Big|\sum_{m\in\cM}\alpha_m K(mn)\Big|^2\\
&\leqslant \frac{\|\boldsymbol\beta\|^2}{|\cN_0||I|}\sum_{n\in\cN_2}\sum_{a\in\cN_0}\sum_{b\in I}\Big|\sum_{m\in\cM}\alpha_m K(m(n+ab))\Big|^2.
\end{align*}
Squaring out we find
\begin{align*}
|\cB(\boldsymbol\alpha,\boldsymbol\beta;K)|^2
&\leqslant \frac{\|\boldsymbol\beta\|^2}{|\cN_0||I|}\mathop{\sum\sum}_{m_1,m_2\in\cM}\sum_{n\in\cN_2}\sum_{a\in\cN_0}\Big|\sum_{b\in I}K(m_1(n+ab))\overline{K(m_2(n+ab))}\Big|.
\end{align*}
The diagonal term with $m_1=m_2$ contributes at most
\begin{align*}
\ll \|\boldsymbol\beta\|^2|\cM||\cN_2|\ll \|\boldsymbol\beta\|^2|\cM||\cN|.
\end{align*}
This gives
\begin{align*}
|\cB(\boldsymbol\alpha,\boldsymbol\beta;K)|^2
&\ll \frac{\|\boldsymbol\beta\|^2}{|\cN_0||I|}\mathop{\sum\sum}_{\substack{m_1,m_2\in\cM\\ m_1\neq m_2}}\sum_{n\in\cN_2}\sum_{a\in\cN_0}\Big|\sum_{b\in I}K(m_1(n+ab))\overline{K(m_2(n+ab))}\Big|\\
&\ \ \ \ +\|\boldsymbol\beta\|^2|\cM||\cN|.
\end{align*}

Put
\begin{align*}
\sigma(x_1,x_2,y)=\mathop{\mathop{\sum\sum}_{m_1\neq m_2\in\cM}\sum_{n\in\cN_2}\sum_{a\in\cN_0}}_{m_1a=x_1,~m_2a=x_2,~n=ay\text{~in }\bF_p}1.
\end{align*}
Hence
\begin{align*}
|\cB(\boldsymbol\alpha,\boldsymbol\beta;K)|^2
&\leqslant \frac{\|\boldsymbol\beta\|^2}{|\cN_0||I|}\mathop{\sum\sum\sum}_{\substack{x_1,x_2,y\in\bF_p\\ x_1\neq x_2}}\sigma(x_1,x_2,y)\Big|\sum_{b\in I}K(x_1(y+b))\overline{K(x_2(y+b))}\Big|\\
&\ \ \ \ +\|\boldsymbol\beta\|^2|\cM||\cN|.
\end{align*}

By H\"older inequality, for all $r\geqslant2$ it follows that
\begin{align}\label{eq:B-P1P2P3}
|\cB(\boldsymbol\alpha,\boldsymbol\beta;K)|^2
&\leqslant \frac{\|\boldsymbol\beta\|^2}{|\cN_0||I|}(\varPi_1)^{1-\frac{1}{r}}(\varPi_2\varPi_3)^{\frac{1}{2r}}+\|\boldsymbol\beta\|^2|\cM||\cN|,
\end{align}
where
\begin{align*}
\varPi_1=\mathop{\sum\sum\sum}_{x_1,x_2,y\in\bF_p}\sigma(x_1,x_2,y),\ \ 
\varPi_2=\mathop{\sum\sum\sum}_{x_1,x_2,y\in\bF_p}\sigma(x_1,x_2,y)^2,
\end{align*}
and
\begin{align*}
\varPi_3&=\mathop{\sum\sum\sum}_{\substack{x_1,x_2,y\in\bF_p\\ x_1\neq x_2}}\Big|\sum_{b\in I}K(x_1(y+b))\overline{K(x_2(y+b))}\Big|^{2r}.
\end{align*}

Trivially we have
\begin{align}\label{eq:Pi1}
\varPi_1\leqslant  |\cM|^2|\cN_0||\cN_2|.
\end{align}
On the other hand, we find
\begin{align*}
\varPi_2
&=\mathop{\sum\cdots\sum}_{\substack{m_1,m_2,m_1',m_2'\in\cM,~n,n'\in\cN_2,~a,a'\in\cN_0\\ m_1a=m_1'a',~m_2a=m_2'a',~na'=n'a\text{~in~}\bF_p}}1\leqslant  |\cM|^2E(\cN_0,\cN_2)\leqslant |\cM|^2\sqrt{E(\cN_0)E(\cN_2)}.
\end{align*}
From Lemma \ref{lm:Shkredov-Kerr} it follows that
\begin{align}\label{eq:Pi2}
\varPi_2
&\ll |\cM|^2\Big(\frac{|\cN_0|^4}{p}+|\cN_0|^{\frac{5}{2}}\Big)^{\frac{1}{2}}
\Big(\frac{|\cN_2|^4}{p}+|\cN_2|^{\frac{5}{2}}\Big)^{\frac{1}{2}}\log^2p.
\end{align}
It remains to estimate $\varPi_3,$ which goes beyond the capacities of Propositions \ref{prop:noncorrelation} and \ref{prop:correlation}.

To bound $\varPi_3,$ we quote a deep result by Kowalski, Michel and Sawin \cite[Theorem 4.3]{KMS20}, from which one may see why we are not able to beat P\'olya--Vinogradov for general trace functions in the sense of Theorem \ref{thm:hyperKl}. We put
\begin{align}\label{eq:Pi(K,bb)}
\varPi(K,\bb):=\mathop{\sum\sum\sum}_{\substack{x_1,x_2,y\in\bF_p\\ x_1\neq x_2}}\prod_{1\leqslant j\leqslant r}K(x_1(y+b_j))\overline{K(x_2(y+b_j))}\overline{K(x_1(y+b_{j+r}))}K(x_2(y+b_{j+r})).
\end{align}

\begin{lemma}\label{lm:KMS-hyperKl}
Let $r\geqslant2$ be an integer.
There exist affine varieties
\begin{align*}
\cV\subseteq \cW\subseteq \bA_{\bZ}^{2r}
\end{align*}
defined over $\bZ$ such that
\begin{align*}
\mathrm{codim}(\cV)=r, \quad \mathrm{codim}(\cW) \geqslant(r-1)/2
\end{align*}
which have the following property: for any $a\in\bF_p^{\times}$ and $k\geqslant2,$ and for all $\bb \in\bF_p^{2r}$ and $K=\kl_k(a\cdot,p),$
we have
\begin{align*}
\varPi(K,\bb)\ll
\begin{cases}
p^3, \ \ & \text { if } \bb\in \cV(\bF_p),\\
p^2,& \text { if } \bb\in (\cW-\cV)(\bF_p),\\
p^{3/2}, & \text { if } \bb\not\in \cW(\bF_p).
\end{cases}
\end{align*}
In all cases, the implied constant depends only polynomially on $k$.
\end{lemma}

The subsequent application of Lemma \ref{lm:KMS-hyperKl} would also require to bound the number of integral points in a box that satisfy a system of polynomial equations in finite fields. See \cite[Lemma 1.7]{Xu20} for details.

\begin{lemma}\label{lm:Xu-counting}
Let $k$ be a positive integer and let $A>0$. Let $X\subseteq \bA_\bZ^k$ be an algebraic variety of dimension $d \geqslant 0$ given by the vanishing of $\leqslant A$ polynomials of degree $\leqslant A$. Let $p$ be a prime number and $0 \leqslant B<p / 2$ an integer. Then
\begin{align*}
|\{\bx\in \bF_p^k: \bx\in X(\bF_p)\cap[B,2B]^k\}|\ll B^d
\end{align*}
where the implied constant depends only on $k$ and $A.$
\end{lemma}

Note that
\begin{align*}
\varPi_3&\leqslant \sum_{\bb\in I^{2r}}|\varPi(K,\bb)|
\end{align*}
with $\varPi(K,\bb)$ given by \eqref{eq:Pi(K,bb)}. According to the location of $\bb,$ we may bound $\varPi(K,\bb)$
appealing to in different situations. According to the existences of $\cV,\cW$ in Lemma \ref{lm:KMS-hyperKl}, we infer from Lemma \ref{lm:Xu-counting} that
\begin{align*}
|\cV(\bF_p)\cap I^{2r}|\ll |I|^{2r-r}=|I|^{r}
\end{align*}
and
\begin{align*}
|(\cW-\cV)(\bF_p)\cap I^{2r}|\ll |I|^{2r-\frac{r-1}{2}}=|I|^{\frac{3r+1}{2}}.
\end{align*}
Hence
\begin{align*}
\varPi_3&\ll |I|^{2r}p^{\frac{3}{2}}+|I|^{\frac{3r+1}{2}}p^2+|I|^rp^3
\end{align*}
for all $r\geqslant2,$
from which and \eqref{eq:Pi1}, \eqref{eq:Pi2} and \eqref{eq:B-P1P2P3}, we infer
\begin{align*}
\cB(\boldsymbol\alpha,\boldsymbol\beta;K)^2
&\ll\frac{\|\boldsymbol\beta\|^2}{|\cN_0||I|} |\cM|^{2-\frac{1}{r}}(|\cN_0||\cN_2|)^{1-\frac{1}{r}}\Big(\frac{|\cN_0|^4}{p}+|\cN_0|^{\frac{5}{2}}\Big)^{\frac{1}{4r}}
\Big(\frac{|\cN_2|^4}{p}+|\cN_2|^{\frac{5}{2}}\Big)^{\frac{1}{4r}}\\
&\ \ \ \times (|I|^{2r}p^{\frac{3}{2}}+|I|^{\frac{3r+1}{2}}p^2+|I|^rp^3)^{\frac{1}{2r}}(\log p)^{\frac{1}{r}}+\|\boldsymbol\beta\|^2|\cM||\cN|.
\end{align*}
The desired estimate then follows by taking $|I|=p^{\frac{3}{2r}}$.

\smallskip

\section{Equidistributions: Proofs of Theorems \ref{thm:Kl-equidistribution} and \ref{thm:ellipticcurve-equidistribution}}
\label{sec:equidistribution-proof}

Given the estimates in Theorem \ref{thm:arbitraryKMN}, the approach to proving equidistributions in Theorems 
\ref{thm:Kl-equidistribution} and \ref{thm:ellipticcurve-equidistribution} is standard.
For the partial completeness, we give the details of Theorem
\ref{thm:Kl-equidistribution} and omit those of the other.

We first quote from \cite[Corollary 3.2]{BSSV18} the following quantitative version of Weyl's criterion for Sato--Tate distributions.

\begin{lemma}\label{lm:BSSVapproximation}
Let $s$ be a positive integer. For any $s$-tuples $(\theta_1, \cdots, \theta_s) \in[0, \pi]^s,$ if there exist some $A,\varDelta>0$ such that
\begin{align*}
\Big|\sum_{1\leqslant i\leqslant s}\sym_k(\theta_i)\Big|\leqslant k^As\varDelta
\end{align*}
holds for every positive integer $k,$ then we have
\begin{align*}
|\{1\leqslant i\leqslant s:2\cos\theta_i\in I\}|=\frac{s}{2\pi}\int_I\sqrt{4-x^2}\ud x+O(s\varDelta^{\frac{1}{A+1}})
\end{align*}
uniformly for the interval $I\subseteq[-2,2].$
\end{lemma}

To prove Theorem \ref{thm:Kl-equidistribution}, we apply Lemma \ref{lm:BSSVapproximation} with
\begin{align*}
s=|\cM||\cN|,\ \ \{\theta_i\}_{1\leqslant i\leqslant s}=\{\theta_p(amn)\}_{(m,n)\in\cM\times\cN},
\end{align*}
and Theorem \ref{thm:arbitraryKMN} guarantees that one may choose 
\begin{align*}
\varDelta=|\cM|^{-\frac{1}{2r}}(p^{\frac{1}{2r}}|\cN|^{-\frac{1}{2}}+p^{\frac{1}{4r}}),
\end{align*}
and $A$ to be an absolute constant. Therefore, we have
\begin{align*}
\frac{1}{|\cM||\cN|}\mathop{\sum\sum}_{\substack{m\in\cM,n\in\cN\\ \kl_2(amn)\in I}}1
&=\frac1{2\pi} \int_I\sqrt{4-x^2}\ud x+O(p^{-\varepsilon^2}),
\end{align*}
provided that
\begin{align*}
|\cM|>p^{\frac{1}{2}+2A\varepsilon},\ \ |\cN|>p^{\varepsilon},\ \ r=[1/\varepsilon].
\end{align*}
Now Theorem \ref{thm:Kl-equidistribution} follows from the arbitrariness of $\varepsilon.$

\smallskip

\bibliographystyle{plainnat}

\end{document}